\documentclass[10pt]{amsart}
\usepackage{amsthm,amsmath,amssymb}

\global\overfullrule10pt

\theoremstyle{plain}
\newtheorem{theorem}{Theorem}
\newtheorem{corollary}[theorem]{Corollary}
\newtheorem{lemma}[theorem]{Lemma}
\newtheorem{proposition}[theorem]{Proposition}
\newtheorem*{theorem*}{Theorem}
\newtheorem*{corollary*}{Corollary}
\newtheorem*{proposition*}{Proposition}
\theoremstyle{definition}

\theoremstyle{remark}
\newtheorem{remark}[theorem]{Remark}

\newcommand\CC{{\mathbf C}}
\newcommand\RR{{\mathbf R}}

\newcommand\NN{{\mathbf N}}
\newcommand\BB{{\mathbf B}}
\newcommand\spr[1]{\langle#1\rangle}

\newcommand\Bn{{\BB^n}}
\newcommand\pBn{{\partial\Bn}}
\newcommand\bn{{B^n}}
\newcommand\pbn{{\partial\bn}}
\renewcommand\dh{\Delta_h}

\newcommand\FF[5]{{}_#1\!F_#2\Big(\begin{matrix}#3\\#4\end{matrix}\Big|#5\Big)}
\newcommand\Hh{$H$-harmonic }
\newcommand\Mh{$M$-harmonic }
\newcommand\chm{\mathcal H^m}
\newcommand\bhm{\mathbf H^m}
\newcommand\On{{O(n)}}
\newcommand\Un{{U(n)}}

\newcommand\HH{{\mathcal H}}

\makeatletter
\newcommand{\vast}{\bBigg@{4}}
\newcommand{\Vast}{\bBigg@{5}}
\makeatother

\renewcommand\[{\begin{equation}}
\renewcommand\]{\end{equation}}
\allowdisplaybreaks[4]

\newcommand\jedna{\mathbf 1}

\newcommand\bG{\mathbf G}
\newcommand\bphi{\varphi}
\newcommand\cZ{\mathcal Z}

\begin{document}

\title[Moebius invariant product]{Uniqueness of $H$-harmonic Moebius invariant inner products on the ball}
\author[P.~Blaschke]{Petr Blaschke}
\address{Mathematics Institute, Silesian University in Opava,
 Na~Rybn\'\i\v cku~1, 74601~Opava, Czech Republic}
 \email{Petr.Blaschke{@}math.slu.cz}
\author[M.~Engli\v s]{Miroslav Engli\v s}
\address{Mathematics Institute, Silesian University in Opava,
 Na~Rybn\'\i\v cku~1, 74601~Opava, Czech Republic {\rm and }
 Mathematics Institute, \v Zitn\' a 25, 11567~Prague~1,
 Czech Republic}
\email{englis{@}math.cas.cz}
\thanks{Research supported by 
 RVO funding for I\v CO 47813059 and~67985840.}
\begin{abstract} We~prove uniqueness of the Moebius invariant semi-inner product
on hyperbolic-harmonic functions on the unit ball of the real n-space,
i.e. on functions annihilated by the hyperbolic Laplacian on the ball.
\end{abstract}

\maketitle

\section{Introduction}
Let $\bn$ be the unit ball in~$\RR^n$, $n>2$. The~orthogonal transformations
\[ x \longmapsto Ux, \qquad x\in\RR^n, \; U\in\On,  \label{TB} \]
map both $\bn$ and its boundary $\pbn$ (the unit sphere) onto themselves,
and so do the Moebius transformations
\[ \phi_a(x) := \frac{a|x-a|^2+(1-|a|^2)(a-x)}{1-2\spr{x,a}+|x|^2|a|^2}, \qquad x\in\bn,  \label{TA} \]
interchanging the origin $0\in\RR^n$ with some point $a\in\bn$.
The~group $G$ generated by the $\phi_a$, $a\in\bn$, and $\On$ via composition is called the Moebius group of~$\bn$.
The~hyperbolic metric
$$ ds = \frac{2|dx|}{1-|x|^2}  $$
on $\bn$ is invariant under the action of~$G$ (and~is the unique Riemannian metric with
this property, up~to a constant factor), and so is the hyperbolic Laplacian
\[ \dh f(x) := (1-|x|^2) [(1-|x|^2)\Delta f(x) + 2(n-2) \spr{x,\nabla f(x)}].  \label{TC}  \]
Functions on $\bn$ annihilated by $\dh$ are called hyperbolic-harmonic, or \Hh for short.
The~vector space $\HH$ of all \Hh functions~is, by~its very definition, invariant under composition
with elements of~$G$, i.e.~with transformations \eqref{TA} and \eqref{TB} above.

In~some sense, \Hh functions are a real-variable counterpart of the \Mh functions on the
unit ball $\Bn$ of~$\CC^n$, whose theory is beautifully presented in Walter Rudin's
classical monograph~\cite{Ru}. This class of functions is preserved, analogously,
by~compositions with elements of the ``complex version'' $\bG$ of the Moebius group,
generated by the unitary maps of $\CC^n$ and by the holomorphic analogues $\bphi_a$
of the maps~$\phi_a$, which again map $\Bn$ onto itself and interchange $a\in\Bn$
with the origin $0\in\CC^n$.

Due to the holomorphy of the elements of~$\bG$,
another class preserved by the compositions just mentioned are the holomorphic functions on~$\Bn$.
A~notable result in the theory of holomorphic functions on the complex unit ball is the fact~that,
up~to a constant factor, there exists a unique semi-inner product on the holomorphic functions
on $\Bn$ which is invariant under the above action of the Moebius group~$\bG$:
namely, the Dirichlet product, given~by
$$ \spr{f,g} = \sum_{\alpha\text{ multi-index}} |\alpha|\frac{\alpha!}{|\alpha|!}f_\alpha\overline{g_\alpha},  $$
where $f(z)=\sum_\alpha f_\alpha z^\alpha$, $g(z)=\sum_\alpha g_\alpha z^\alpha$ are the Taylor expansions
of $f$ and~$g$, respectively, and we are using the usual multi-index notation.

It~has been an open question --- mentioned explicitly also on p.~180 of Manfred Stoll's monograph~\cite{St} ---
whether an analogous result also holds in the real-variable setting, i.e.~whether there exists, similarly,
a~Moebius-invariant semi-inner product on the space of \Hh functions on~$\bn$, $n>2$, and whether it is unique.

The existence of such ``hyperbolic Dirichlet'' semi-inner product has recently been shown
by the current authors and E.-H.~Youssfi in~\cite{BEY}, by~a~kind of ``analytic continuation''
of the corresponding weighted Bergman inner products, in the spirit of Rossi and~Vergne.
In~this paper, we~complement this construction by proving that, up~to a constant factor,
this $G$-invariant semi-inner product on \Hh functions is again unique.
That~is, the situation for \Hh functions on the real ball $\bn$ is completely parallel
to that for the holomorphic functions on the complex ball~$\Bn$.

In~the holomorphic case, the proof of uniqueness relies heavily on manipulations with Taylor coefficients,
see e.g.~Section~6.4 in the book of Kehe Zhu~\cite{Zh}. This approach seems not to be directly applicable
for the real ball; instead, our~method proceeds by evaluation of certain integrals involving zonal harmonics.
These may be of independent interest and may find potential applications elsewhere; besides, the~same idea
can also be used to give another proof of uniqueness also in the holomorphic case.

In~Section~\ref{dva} below, we~collect some preliminaries on \Hh functions and zonal harmonics.
Section~\ref{tri} contains technical computations involving the zonal polynomials, which are then
used in Section~\ref{ctyri} to give the proof of uniqueness of the ``hyperbolic-Dirichlet'' semi-inner product.
In~the final section, Section~\ref{pet}, we~briefly present the alternative proof for the holomorphic case.

\section{Notation and preliminaries} \label{dva}
Each orthogonal transformation $U\in\On$ maps the unit sphere $\pbn$ onto itself,
and the normalized surface measure $d\sigma$ on $\pbn$ is invariant under~$U$.
It~follows that the composition with elements of $\On$,
\[ T_U : f\longmapsto f\circ U^{-1}  \label{TD}  \]
is a unitary representation of $\On$ on~$L^2(\pbn,d\sigma)$. The~decomposition of this
representation into irreducible subspaces is realized by the spherical harmonics;
see e.g.~\cite{ABR}, especially Chapter~5. Namely, for~an integer $m\ge0$,
let $\chm$ be the space of restrictions to the unit sphere $\pbn$ of harmonic polynomials
on~$\RR^n$ homogeneous of degree~$m$. Then $\chm$ is invariant under the action~\eqref{TD},
is~$\On$-irreducible (i.e.~has no proper $\On$-invariant subspaces) and
\[ L^2(\pbn,d\sigma) = \bigoplus_{m=0}^\infty \chm .  \label{TE}  \]

Associated to each $\chm$ is the space of ``solid harmonics''
$$ \bhm = \{f\in C(\overline{\bn}): \;f\text{ is \Hh on }\bn \text{ and }f|_\pbn\in\chm\} .  $$
It~is then a fact that any $f_m\in\bhm$ is necessarily of the form
\[ f_m(r\zeta) = S_m(r^2) r^m f_m(\zeta), \quad f_m|\pbn\in\chm, \qquad 0\le r\le1, \; \zeta\in\pbn,  \label{TF} \]
where
\[ S_m(t) := \frac{(n-1)_m}{(\frac n2)_m} \FF21{m,1-\tfrac n2}{m+\frac n2}t  \label{TG}  \]
with the Gauss hypergeometric function $_2F_1$ and the Pochhammer symbol $(x)_m:=x(x+1)\dots(x+m-1)$.
Furthermore, any \Hh function $f$ on $\bn$ can be uniquely written in the form
\[ f = \sum_m f_m, \qquad f_m\in\bhm,   \label{TH}  \]
with the sum converging uniformly on compact subsets of~$\bn$.
The~components $f_m$ can be recovered from $f$ by
\[ f_m(r\zeta) = \int_\pbn f(r\eta) Z_m(\zeta,\eta) \,d\sigma(\eta),  \label{YH}  \]
where the zonal harmonics $Z_m(x,y)$ are defined by the generating function
$$ \sum_{m=0}^\infty t^m Z_m(x,y) = \frac{1-t^2|x|^2|y|^2}{(1-2t\spr{x,y}+t^2|x|^2|y|^2)^{n/2}}, \qquad \|x\|,\|y\|\le1, \;0\le t<1.  $$
In~particular, $Z_m(x,y)$ is homogeneous of degree $m$ in both $x$ and~$y$, and has the reproducing property
\[ f_m(x) = \int_\pbn f_m(y) Z_m(x,y) \,d\sigma(y), \qquad x\in\pbn, \; f_m\in\chm.  \label{TI}  \]
See Chapter~8 in \cite{ABR} for a nice exposition on spherical harmonics.

The reader is referred to \cite{St} for a detailed exposition of the facts above,
with a convenient overview and further developments in \cite{St2} and~\cite{Ur}.
Note however that, for ease of notation, our $S_m(t)$ is $S_m(\sqrt t)$ in the
notation of these references.

\section{Preparatory calculations} \label{tri}
We~proceed to compute an integral involving zonal harmonics which we will need later.

\begin{lemma} \label{mscalarlemma}
For $x,y\in \pbn $, i.e. $|x|=|y|=1$, we~have
\begin{align*}
\spr{\phi_a(x),\phi_a(y)}=1-\frac{(1-\spr{x,y})(1-|a|^2)^2}{(1-2\spr{a,x}+|a|^2)(1-2\spr{a,y}+|a|^2)}.
\end{align*}
\end{lemma}

\begin{proof} This follows by a tedious yet totally routine calculation. \end{proof}

\begin{lemma} \label{Zonalthreereq}
For $x,y\in\RR^n$ such that $|x|=|y|=1$,
$$ 2\spr{x,y} Z_m(x,y)=\frac{m+1}{m+\frac n2}Z_{m+1}(x,y)+\frac{m+n-3}{m+\frac n2-2}Z_{m-1}(x,y). $$
\end{lemma}

\begin{proof}
The well known three term recurrence relation for the underlying ultraspherical (or Gegenbauer) polynomials
cf.~\cite[Theorem~5.38]{ABR} and \cite[\S10.9~(18) and~(13)]{BE} says that
$$ 2\spr{x,y} Z_m(x,y)=\frac{m+1}{m+\frac n2}Z_{m+1}(x,y)+\frac{m+n-3}{m+\frac n2-2}|x|^2|y|^2 Z_{m-1}(x,y). $$
Specializing to $|x|=|y|=1$ yields the lemma.
\end{proof}

\begin{lemma} \label{Le3}
For any $k\in \NN$,
\begin{align*}
Z_m(y,\nabla_u)|u|^{2k}&=m!\binom km 2^m |u|^{2(k-m)}Z_m(y,u).
\end{align*}
\end{lemma}

\begin{proof} Using in turn the Leibnitz rule, the fact that $Z_m(x,y)$ is harmonic in each argument, and the binomial theorem, we~have
\begin{align*}
Z_m(y,\nabla_u)|u|^{2k}&=Z_m(y,\nabla_u)\Delta_x^k e^{\spr{u,x}}|_{x=0}\\
&=\Delta_x^k Z_m(y,x)e^{\spr{u,x}}|_{x=0}\\
&=(|u|^2+2\spr{u,\nabla_x} +\Delta_x)^k Z_m(y,x)|_{x=0}\\
&=(|u|^2+2\spr{u,\nabla_x})^k Z_m(y,x)|_{x=0}\\
&=\sum_{j=0}^k\binom kj |u|^{2(k-j)} 2^j \spr{u,\nabla_x}^j Z_m(y,x)|_{x=0}.
\end{align*}
Now, since $Z_m$ is homogeneous polynomial of degree $m$ we have
$$
\spr{u,\nabla_x}^j Z_m(y,x)|_{x=0}=\delta_{j,m} m! Z_m(y,u).
$$
Thus
$$ Z_m(y,\nabla_u)|u|^{2k} = \binom km |u|^{2(k-m)} 2^m m!Z_m(y,u), $$
as claimed.
\end{proof}

\begin{corollary} \label{Zme}
Let $u,y\in\RR^n$. Then
$$ \int_\pbn{} Z_m(x,y)e^{\spr{u, x}}\,d\sigma(x) =
 \frac{Z_m(u,y)}{(\frac n2)_m 2^m} \FF01{-}{\frac n2+m}{\frac{|u|^2}4}.  $$
\end{corollary}

\begin{proof} Expanding the exponential into power series and integrating term by term shows that
$$ \int_\pbn e^{\spr{u,x}}\,d\sigma(x) = \FF01{-}{\frac n2}{\frac{|u|^2}4} = \sum_{k=0}^\infty \frac{|u|^{2k}}{(\frac n2)_k 4^k k!}. $$
Differentiating under the integral sign, we~get by Lemma~\ref{Le3}
\begin{align*}
& \int_\pbn{} Z_m(y,x) e^{\spr{u,x}}\,d\sigma(x) = \int_\pbn Z_m(y,\nabla_u) e^{\spr{u,x}}\,d\sigma(x) \\
&\hskip4em =Z_m(y,\nabla_u)\int_\pbn e^{\spr{u,x}}\,d\sigma(x) = Z_m(y,\nabla_u) \FF01{-}{\frac n2}{\frac{|u|^2}4} \\
&\hskip4em =\sum_{k=0}^\infty \frac{Z_m(y,\nabla_u)|u|^{2k}}{(\frac n2)_k 4^k k!}
 =Z_m(y,u)\sum_{k=m}^\infty \frac{|u|^{2(k-m)}2^m}{(\frac n2)_k 4^k (k-m)!} \\
&\hskip4em = Z_m(y,u)\sum_{k=0}^\infty \frac{|u|^{2k}}{(\frac n2)_{k+m}4^{k+m} k!}
 =\frac{Z_m(y,u)}{2^m(\frac n2)_m}\sum_{k=0}^\infty \frac{|u|^{2k}}{(\frac n2+m)_k 4^k k!} \\
&\hskip4em = \frac{Z_m(y,u)}{2^m(\frac n2)_m} \FF01{-}{\frac n2+m}{\frac{|u|^{2}}4},
\end{align*}
as claimed.
\end{proof}

\begin{corollary} \label{Zmee}
Let $u,v\in\RR^n$. Then
$$ \int_\pbn \int_\pbn Z_m(x,y)e^{\spr{u,x}}e^{\spr{v,y}}\,d\sigma(x)\,d\sigma(y) =
\frac{Z_m(u,v)}{(\frac n2)_m^2 4^m} \FF01{-}{\frac n2+m}{\frac{|u|^2}4} \FF01{-}{\frac n2+m}{\frac{|v|^2}4}.  $$
\end{corollary}

\begin{proof} Apply Lemma~\ref{Zme} twice. \end{proof}

\begin{lemma} \label{hyp1} For $b,c>0$ and $-1<x<1$, we have
$$ \int_0^\infty t^{b-1} \FF01{-}{c}{\frac{t^2 x}{(1+x)^2}} e^{-t}\,dt=
\Gamma(b)(1+x)^b \FF21{b,1+b-c}{c}x.  $$
\end{lemma}

\begin{proof}
Expanding the ${}_0\!F_1$ function into its Taylor series and integrating term by term yields
$$ \int_0^\infty t^{b-1} \FF01{-}{c}{\frac{t^2 x}{(1+x)^2}} e^{-t}\,dt =
 \Gamma(b) \FF21{\frac b2,\frac{b+1}2}{c}{\frac{4x}{(1+x)^2}}.  $$
Now just apply the well known transform for ${}_2\!F_1$ functions \cite[15.8.15]{NIST}.
\end{proof}

\begin{corollary}\label{Cor6}
For $a\in\bn$,
$$ \int_\pbn \int_\pbn \frac{Z_m(x,y)\,d\sigma(x)\,d\sigma(y)}{(1+|a|^2-2\spr{a,x})(1+|a|^2-2\spr{a,y})}
 =\frac{(n-2)_m m!}{(\frac n2)_m(\frac n2-1)_m }|a|^{2m} \FF21{m+1,2-\frac n2}{\frac n2+m}{|a|^2}^2.  $$
\end{corollary}

\begin{proof} The proof stems from the following representation
$$ \frac1{(1+|a|^2-2\spr{a,x})(1+|a|^2-2\spr{a,y})}
 =\frac1{(1+|a|^2)^2} \int_0^\infty \int_0^\infty e^{\alpha \spr{u,x}} e^{\beta \spr{u,y}} e^{-\alpha} e^{-\beta}\,d\alpha\,d\beta , $$
where
$$ u=\frac{2a}{1+|a|^2}. $$
Interchanging the order of integration and applying Corollary~\ref{Zmee}, we~get
\begin{align*}
&\int_\pbn \int_\pbn \frac{Z_m(x,y)}{(1+|a|^2-2\spr{a,x})(1+|a|^2-2\spr{a,y})} \,d\sigma(x)\,d\sigma(y)\\
&= \frac1{(1+|a|^2)^2} \int_0^\infty \int_0^\infty
 \frac{(\alpha\beta)^m Z_m(u,u)}{(\frac n2)_m^2 4^m} \FF01{-}{\frac n2+m}{\frac{\alpha^2|u|^2}4} \FF01{-}{\frac n2+m}{\frac{\beta^2|u|^2}4}
 e^{-\alpha}e^{-\beta}\,d\alpha\,d\beta\\
&= \frac{Z_m(u,u)m!^2}{(\frac n2)_m^2 4^m}(1+|a|^2)^{2m} \FF21{m+1,2-\frac n2}{\frac n2+m}{|a|^2}^2.
\end{align*}
where in the last equality we have used Lemma~\ref{hyp1}. Now it suffices to realize that
$$ Z_m(u,u)=|u|^{2m}\frac{(\frac n2)_m(n-2)_m}{(\frac n2-1)_m m!}=\frac{4^m |a|^{2m}}{(1+|a|^2)^{2m}}\frac{(\frac n2)_m(n-2)_m}{(\frac n2-1)_m m!}. $$
\end{proof}

\begin{proposition} \label{Puniq} For any $m\in\NN$ and $a\in\bn$,
\begin{align*}
&\int_\pbn \int_\pbn Z_m(x,y)\spr{\phi_a(x),\phi_a(y)}\,d\sigma(x)\,d\sigma(y)\\
&=\delta_{m,0}-(1-|a|^2)^{2n-4}\Big[\Phi_m(|a|^2)-\frac{m+1}{2m+n}\Phi_{m+1}(|a|^2)-\frac{m+n-3}{2m+n-4}\Phi_{m-1}(|a|^2)\Big],
\end{align*}
where
$$ \Phi_m(x):=\frac{(n-2)_m m!}{(\frac n2)_m(\frac n2-1)_m}x^m \FF21{\frac n2-1,m+n-2}{m+\frac n2}x ^2,\qquad \Phi_{-1}(x):=0. $$
\end{proposition}

\begin{proof} Combine Lemma~\ref{mscalarlemma}, Lemma~\ref{Zonalthreereq}
and Corollary~\ref{Cor6} and apply the~Euler transform (cf.~\cite[\S2.1~(22)]{BE})
$$ \FF21{a,b}cz = (1-z)^{c-a-b}\FF21{c-a,c-b}cz  $$
to the resulting ${}_2\!F_1$ functions.
\end{proof}

\section{Uniqueness} \label{ctyri}
Let~us now look for semi-inner products (and the corresponding \hbox{semi-}norms) on spaces of
\Hh functions on the unit ball $\bn$ of the real $n$-space~$\RR^n$, $n\ge2$, which are invariant
under Moebius maps $\phi_a$ from \eqref{TA} as~well as under orthogonal rotations~\eqref{TB}.
By~Schur's lemma from representation theory (see e.g.~Proposition~1 in \cite{EY} for the detailed
argument in the setting of the unitary group $U(n)$ instead of~$O(n)$), the~rotation invariance
implies that any such semi-norm which is finite on each $\bhm$ must come in the form
\[ \|f\|^2 = \sum_m c_m \|f_m\|^2_\pbn,  \label{YA}  \]
where $f(z)=\sum_m f_m(z)$ is the decomposition \eqref{TH} of an \Hh function $f$ into its Peter-Weyl components
under the action of~$\On$ recalled in the preceding section, $\|\cdot\|_\pbn$ denotes the norm in $L^2(\pbn,d\sigma)$,
and $c_m\ge0$ are some coefficients. For~typographic reasons, we~often write $[f]_m$ instead of~$f_m$.
By~polarization, the~associated semi-inner product is given~by
$$ \spr{f,g} = \sum_m c_m \spr{f_m,g_m}_\pbn  $$
for any $f,g$ whose semi-norm \eqref{YA} is finite.

It~was shown in \cite{BEY} that the choice $c_m=(n-1)_m/\Gamma(m)$, that~is,
\[ \|f\|_\circ^2 := \sum_{m>0} \frac{(n-1)_m}{\Gamma(m)} \|f_m\|^2_\pbn ,  \label{YM}  \]
leads to a semi-norm which is Moebius invariant, i.e.~$\|f\circ\phi\|_\circ=\|f\|_\circ$ for any $\phi\in G$.
Our~main result in this paper is the following.

\begin{theorem} If~a semi-norm of the form \eqref{YA} is Moebius invariant, then it is a constant multiple of~$\|\cdot\|_\circ$.
\end{theorem}

\begin{proof} Invariance under Moebius maps $\phi_a$ means that
$$ \|f\circ\phi_a\| = \|f\|   \qquad \forall a\in\bn  $$
for all \Hh~$f$ (both sides being allowed to be possibly infinite). Since the functions~$z_j S_1(|z|^2)$,
$j=1,\dots,n$, are~\Hh by~\eqref{TF}, this in particular implies
\[ \|\phi_{a,j} S_1(|\phi_a|^2)\| \qquad\text{is independent of }a\in\bn,  \label{YG} \]
where $\phi_{a,j}$ denotes the $j$-th coordinate of~$\phi_a$. Summing over $j$ and using~\eqref{YA},
it~follows that the function
\[ F(a) := \sum_{j=1}^n \sum_m c_m \|[\phi_{a,j}]_m\|^2_\pbn , \qquad a\in\bn,  \label{YB}  \]
is constant. (Note that $\phi_a$ maps $\pbn$ into itself and $S_1(1)=1$.)
Quite generally, for any $f$ \Hh on~$\bn$ and continuous on the closure, we~have by \eqref{TF} and~\eqref{TI},
\[ f_m(z) = S_m(|z|^2) \int_\pbn f(\zeta) Z_m(z,\zeta) \, d\sigma(\zeta) , \qquad \forall z\in\bn,  \label{YE}  \]
with $Z_m(x,y)$ the zonal harmonic of degree~$m$ from the preceding section.
Feeding this into \eqref{YB} yields
$$ F(a) = \sum_{j=1}^n \sum_m c_m \int_\pbn \int_\pbn \int_\pbn \phi_{a,j}(\zeta) Z_m(\xi,\zeta)
 \overline{\phi_{a,j}(\eta) Z_m(\xi,\eta)} \, d\sigma(\zeta) \,d\sigma(\eta) \,d\sigma(\xi) ,  $$
or, using the reproducing property \eqref{TI} of $Z_m$ again for the $\xi$-integration,
\begin{align*}
F(a) &= \sum_{j=1}^n \sum_m c_m \int_\pbn \int_\pbn \phi_{a,j}(\zeta) Z_m(\eta,\zeta)
 \overline{\phi_{a,j}(\eta)} \, d\sigma(\zeta) \,d\sigma(\eta)  \\
&= \sum_m c_m \int_\pbn \int_\pbn \spr{\phi_a(\zeta),\phi_a(\eta)} Z_m(\eta,\zeta)
 \, d\sigma(\zeta) \,d\sigma(\eta) .
\end{align*}
Using Proposition~\ref{Puniq}, this yields
\begin{align*}
\frac{c_0-F(a)}{(1-|a|^2)^{2n-4}} &= \sum_{m=0}^\infty c_m \Big[\Phi_m(|a|^2)-\frac{m+1}{2m+n}\Phi_{m+1}(|a|^2)-\frac{m+n-3}{2m+n-4}\Phi_{m-1}(|a|^2)\Big] \\
&= \sum_{m=0}^\infty \Big[ c_m - \frac{m}{2m+n-2} c_{m-1} - \frac{m+n-2}{2m+n-2} c_{m+1} \Big] \Phi_m(|a|^2) \\
&= \sum_{m=0}^\infty \Big[ \frac{m}{2m+n-2} (c_m-c_{m-1}) - \frac{m+n-2}{2m+n-2} (c_{m+1}-c_m) \Big] \Phi_m(|a|^2)
\end{align*}
with $c_{-1}:=0$. Setting $a=0$, this becomes $c_0-F(0)=c_0-c_1$, so $F$ being constant means that $F(a)\equiv c_1$.
Introducing the notation $d_m:=c_m-c_{m-1}$, we~can thus continue with
\[ \frac{d_1}{(1-|a|^2)^{2n-4}} = \sum_{m=0}^\infty \Big[\frac{m+n-2}{2m+n-2} d_{m+1} - \frac{m}{2m+n-2} d_m\Big] \Phi_m(|a|^2).  \label{YD}  \]
Denote by $a_{m,k}$ the Taylor coefficients of the function $\Phi_m(x)/x^m$:
$$ \frac{(n-2)_m m!}{(\frac n2)_m(\frac n2-1)_m} \FF21{\frac n2-1,m+n-2}{m+\frac n2}x^2 =: \sum_{k=0}^\infty a_{m,k} x^k . $$
Note that $a_{m,k}>0$, since $n>2$ and all Taylor coefficients of the ${}_2\!F_1$ are therefore strictly positive real numbers.
Comparing the Taylor coefficients at like powers of~$|a|^2$ in~\eqref{YD}, we~obtain for all $k=0,1,2,\dots$,
$$ \frac{(2n-4)_k}{k!}d_1 = \sum_{m=0}^k a_{m,k-m} \Big[\frac{m+n-2}{2m+n-2} d_{m+1} - \frac{m}{2m+n-2} d_m\Big] , $$
or
\begin{align*}
a_{k,0} \frac{k+n-2}{2k+n-2} d_{k+1} &= \frac{(2n-4)_k}{k!}d_1 + a_{k,0} \frac{k}{2k+n-2}d_k \\
&\hskip2em - \sum_{m=0}^{k-1} a_{m,k-m} \Big[\frac{m+n-2}{2m+n-2} d_{m+1} - \frac{m}{2m+n-2} d_m\Big].
\end{align*}
For $k=0$, this reads $a_{0,0}d_1=d_1$, which is always fulfilled since $a_{0,0}=1$. For $k>0$, this gives a recurrence
formula which uniquely determines $d_2,d_3,\dots$ from~$d_1$. Since we know from \eqref{YM} that
\[ d_m = \frac{(n-1)_m}{\Gamma(m)} - \frac{(n-1)_{m-1}}{\Gamma(m-1)} = \frac{n-1}{n-2} \frac{(n-2)_m}{\Gamma(m)}  \label{YI} \]
solves this recurrence for $d_1=n-1$, it~follows that, necessarily,
$$ d_m = \alpha \frac{n-1}{n-2} \frac{(n-2)_m}{\Gamma(m)}, \qquad \forall m\in\NN,  $$
with some constant $\alpha>0$ (explicitly, $\alpha=d_1/(n-1)$), and, consequently,
$$ c_m = c_0 + \alpha \frac{(n-1)_m}{\Gamma(m)} .  $$

On~the other hand, similarly to~\eqref{YG}, we~must have for all $a\in\bn$ and $j=1,\dots,n$,
$$ \spr{\jedna,z_j S_1(|z|^2)} = \spr{\jedna\circ\phi_a,(z_j S_1(|z|^2))\circ\phi_a} = \spr{\jedna,\phi_{a,j}S_1(|\phi_a|^2)},  $$
or, using~\eqref{YA},
$$ \sum_m c_m \spr{[\jedna]_m,[z_j]_m}_\pbn = \sum_m c_m \spr{[\jedna]_m,[\phi_{a,j}]_m}_\pbn , $$
that~is, since $[\jedna]_m=\delta_{m0}\jedna$,
$$ 0 = c_0 \spr{\jedna,[\phi_{a,j}]_0}_\pbn, $$
or, since, quite generally, $[f]_0=f(0)\jedna$ (by~\eqref{YE} and~\eqref{TH}),
$$ 0 = c_0 \phi_{a,j}(0) = c_0 a_j.  $$
As~$a\in\bn$ is arbitrary, necessarily $c_0=0$.

Consequently,
$$ c_m = \alpha \frac{(n-1)_m}{\Gamma(m)},   $$
proving the claim.
\end{proof}

\begin{remark} It~would be amusing to have a direct proof that the quantities \eqref{YI} are a solution to~\eqref{YD}.
This is tantamount~to
$$ (1-z)^{4-2n} = \sum_{m=0}^\infty \frac{(n-2)_m}{m!} \Phi_m(z)   $$
for $z\in\CC$, $|z|<1$.
\end{remark}

\section{The holomorphic case} \label{pet}
In~the holomorphic situation, we~are looking for semi-inner products (and corresponding semi-norms)
on spaces of holomorphic functions on the unit ball $\Bn$ of the complex $n$-space~$\CC^n$, $n\ge1$,
which are invariant under the holomorphic Moebius maps $\bphi_a$ interchanging the origin with $a\in\Bn$,
as~well as under unitary rotations $z\mapsto Uz$, $U\in\Un$. By~Schur's lemma again, the~rotation invariance
implies that any such norm must come in the form
\[ \|f\|^2 = \sum_m c_m \|f_m\|^2_\pBn,  \label{XA}  \]
where $f(z)=\sum_m f_m(z)$ is the decomposition of a holomorphic function $f$ into its homogeneous
components of the Taylor expansion, $\|\cdot\|_\pBn$ denotes the norm in $L^2(\pBn,d\sigma)$ with respect
to the normalized surface measure $d\sigma$ on the sphere~$\pBn$, and $c_m\ge0$ are some coefficients.
For~typographic reasons, we~again often write $[f]_m$ instead of~$f_m$.

Invariance under Moebius maps $\bphi_a$ means that
$$ \|f\circ\bphi_a\| = \|f\|   \qquad \forall a\in\Bn $$
for all holomorphic~$f$ (both sides being allowed to be possibly infinite). Since the coordinate
functions~$z_j$, $j=1,\dots,n$, are~holomorphic, this in particular implies
\[ \|\bphi_{a,j}\| \qquad\text{is independent of }a\in\Bn,  \label{XG} \]
where $\bphi_{a,j}$ denotes the $j$-th coordinate of~$\bphi_a$. Summing over $j$ and using~\eqref{XA},
it~follows that the function
\[ F(a) := \sum_{j=1}^n \sum_m c_m \|[\bphi_{a,j}]_m\|^2_\pBn , \qquad a\in\Bn,  \label{XB}  \]
is constant. Quite generally, for any $f$ holomorphic on~$\Bn$ and continuous on the closure,
\[ f_m(z) = \int_\pBn f(\zeta) \cZ_m(z,\zeta) \, d\sigma(\zeta) ,  \label{XE}  \]
where
\[ \cZ_m(x,y) = \frac{(n)_m}{m!} \spr{x,y}^m   \label{XC}  \]
is the reproducing kernel of the space of homogeneous polynomials of degree $m$ with respect to
the $L^2(\pBn,d\sigma)$ inner product. Feeding this into \eqref{XB} yields
$$ F(a) = \sum_{j=1}^n \sum_m c_m \int_\pBn \int_\pBn \int_\pBn \bphi_{a,j}(\zeta) \cZ_m(\xi,\zeta)
 \overline{\bphi_{a,j}(\eta) \cZ_m(\xi,\eta)} \, d\sigma(\zeta) \,d\sigma(\eta) \,d\sigma(\xi) ,  $$
or, using the reproducing property of $\cZ_m$ for the $\xi$-integration,
\begin{align}
F(a) &= \sum_{j=1}^n \sum_m c_m \int_\pBn \int_\pBn \bphi_{a,j}(\zeta) \cZ_m(\eta,\zeta)
 \overline{\bphi_{a,j}(\eta)} \, d\sigma(\zeta) \,d\sigma(\eta)  \nonumber \\
&= \sum_m c_m \int_\pBn \int_\pBn \spr{\bphi_a(\zeta),\bphi_a(\eta)} \cZ_m(\eta,\zeta)
 \, d\sigma(\zeta) \,d\sigma(\eta) .   \label{XD}
\end{align}
Recall that
\[ 1-\spr{\bphi_a(\zeta),\bphi_a(\eta)} = \frac{(1-|a|^2)(1-\spr{\zeta,\eta})} {(1-\spr{\zeta,a})(1-\spr{a,\eta})}.  \label{XH} \]
Hence
$$ F(a) = \sum_m c_m \Big[ \delta_{m0} - (1-|a|^2) \int_\pBn \int_\pBn
 \frac{(1-\spr{\zeta,\eta}) \cZ_m(\eta,\zeta)} {(1-\spr{\zeta,a})(1-\spr{a,\eta})}
 \, d\sigma(\zeta) \,d\sigma(\eta) \Big]  . $$
Applying \eqref{XE} to the holomorphic function
$$ f(z) = \frac{1-\spr{z,\eta}} {1-\spr{z,a}} = \sum_k [\spr{z,a}^k - (1-\delta_{k0})\spr{z,\eta}\spr{z,a}^{k-1}] ,  $$
the $\zeta$-integration yields
\begin{align*}
F(a) &= \sum_m c_m \Big[ \delta_{m0} - (1-|a|^2) \int_\pBn \frac{\spr{\eta,a}^m -
 (1-\delta_{m0})\spr{\eta,\eta}\spr{\eta,a}^{m-1}} {(1-\spr{a,\eta})} \,d\sigma(\eta) \Big] \\
&= \sum_m c_m \Big[ \delta_{m0} - (1-|a|^2) \int_\pBn \frac{\frac{m!}{(n)_m} \cZ_m(\eta,a) -
 (1-\delta_{m0}) \frac{(m-1)!}{(n)_{m-1}} Z_{m-1}(\eta,a)} {(1-\spr{a,\eta})} \,d\sigma(\eta) \Big] ,
\end{align*}
so applying \eqref{XE} one more time finally gives
\[ F(a) = \sum_m c_m \Big[ \delta_{m0} - (1-|a|^2)
 \Big( \frac{m!}{(n)_m} |a|^{2m} - (1-\delta_{m0}) \frac{(m-1)!}{(n)_{m-1}} |a|^{2m-2} \Big) \Big] .  \label{XF} \]
For $a=0$, this becomes
$$ F(0) = \sum_m c_m [\delta_{m0} - (\delta_{m0} - \delta_{m1})] = c_1 . $$
Subtracting $c_0$ from both sides of $F(a)-F(0)=0$, we~thus get from \eqref{XF}
$$ \frac{c_0-c_1}{1-|a|^2} = \sum_m c_m \Big( \frac{m!}{(n)_m} |a|^{2m} - (1-\delta_{m0}) \frac{(m-1)!}{(n)_{m-1}} |a|^{2m-2} \Big) .  $$
Equating coefficients at like powers of $|a|^2$ on both sides, we~see that
$$ c_0-c_1 = \frac{m!}{(n)_m} (c_m-c_{m+1}) \qquad\forall m.  $$
Letting
$$ b_m := \frac{c_1-c_0}n m \frac{(n)_m}{m!},  $$
we~have
$$ b_m-b_{m+1} = \frac{c_1-c_0}n \frac{(n)_m}{m!} \Big[m-(m+1)\frac{n+m}{m+1} \Big] = (c_0-c_1)\frac{(n)_m}{m!} .  $$
It~follows that $b_m-c_m$ is independent of~$m$, so
$$ c_m = b_m + (c_0-b_0) = c_0 + \frac{c_1-c_0}n m \frac{(n)_m}{m!} .  $$

Similarly as in~\eqref{XG}, we~must have for all $a\in\Bn$ and $j=1,\dots,n$,
$$ \spr{\jedna,z_j} = \spr{\jedna\circ\bphi_a,z_j\circ\bphi_a} = \spr{\jedna,\bphi_{a,j}},  $$
or, using~\eqref{XA},
$$ \sum_m c_m \spr{[\jedna]_m,[z_j]_m}_\pBn = \sum_m c_m \spr{[\jedna]_m,[\bphi_{a,j}]_m}_\pBn , $$
that~is, as $[\jedna]_m=\delta_{m0}\jedna$,
$$ 0 = c_0 \spr{\jedna,[\bphi_{a,j}]_0}_\pBn, $$
or, since, quite generally, $[f]_0=f(0)\jedna$,
$$ 0 = c_0 \bphi_{a,j}(0) = c_0 a_j.  $$
As~$a\in\Bn$ is arbitrary, necessarily $c_0=0$.

Consequently,
$$ c_m = \frac{c_1}n m \frac{(n)_m}{m!} ,   $$
i.e.~up to a constant factor the unique Moebius-invariant inner product on the space of
holomorphic functions is the usual Dirichlet inner product on~$\Bn$.

\end{document}